\documentclass[12pt]{article}
\usepackage{epic,eepic,epsf,epsfig}
\usepackage{amsmath}
\usepackage{amssymb}
\usepackage{amsbsy}

\usepackage{amsfonts}
\newtheorem{theorem}{Theorem}[section]%
\newtheorem{lemma}[theorem]{Lemma}%
\newtheorem{cor}[theorem]{Corollary}%

\def\f{\noindent}

\newcommand{\qed}{\mbox{\raisebox{0.7ex}{\fbox{}}} \vspace{4truemm}}
\def\demo{\f {\bf Proof.}\hskip10pt}

\begin{document}

\baselineskip 16pt

\title{ \vspace{-1.2cm}
A note on an extension of Shidov's theorem and Shlyk's theorem
\thanks{\scriptsize This work was supported in part by NSFC (12571024) and Shandong Provincial Natural Science Foundation, China (ZR2017MA022).
\newline
 \hspace*{0.5cm} \scriptsize $^{\ast\ast}$Corresponding
  author.\newline
 \hspace*{0.5cm} \scriptsize{E-mail addresses:} jtshi2005@pku.org.cn\,(J. Shi),\,\,lifh2026@s.ytu.edu.cn\,(F. Li).}}

\author{Jiangtao Shi\,$^{\ast\ast}$,\,\,Fanghua Li\\
\\
{\small{\em School of Mathematics and Information Sciences,
Yantai University, Yantai 264005, P.R. China}}}

\date{ }

\maketitle \vspace{-.8cm}

\begin{abstract}
Let $G$ be a non-solvable group and $p$ any fixed prime divisor of $|G|$. We prove that the intersection of all non-nilpotent maximal subgroups of $G$
of order divisible by $p$ containing the normalizer of some Sylow subgroup is $p$-closed, which provides a further generalization of Shidov's theorem and Shlyk's theorem.

\medskip
\f {\bf Keywords:} non-solvable group; non-nilpotent maximal subgroup; Sylow subgroup; $p$-closed\\
{\bf MSC(2020):} 20D10; 20D20; 20E28
\end{abstract}

\section{Introduction}

In this paper all groups are assumed to be finite. Recall that the intersection of all maximal subgroups of a group is called a Frattini subgroup. As
generalizations, Shidov {\rm\cite{shidov}} showed that the intersection of all non-nilpotent maximal subgroups of a non-solvable group is nilpotent, which follows easily that
the intersection of all non-nilpotent maximal subgroups of a non-solvable group is equal to its Frattini subgroup. In {\rm\cite{shlyk},
Shlyk proved that the intersection of all non-normal non-nilpotent maximal subgroups of a non-solvable group is nilpotent. Recently, Shi, Xu and Li
{\rm\cite[Theorem 1.2]{shixuli2026}} gave the following result.

\begin{theorem}{\rm\cite[Theorem 1.2]{shixuli2026}}\ \ \label{thx} The intersection of all non-nilpotent maximal subgroups of a non-solvable $G$ containing the normalizer
of some Sylow subgroup is nilpotent.
\end{theorem}

Let $G$ be a group and $p$ any fixed prime divisor of $|G|$. In {\rm\cite[Theorem 1.1]{lishi2026}} Li and Shi provided a complete characterization of a group $G$ in which
every maximal subgroup of order divisible by $p$ containing the normalizer of some Sylow subgroup is nilpotent. In particular, {\rm\cite[Corollary 1.2]{lishi2026}} indicated
that any non-solvable $G$ has a non-nilpotent maximal subgroup of order divisible by $p$ containing the normalizer of some Sylow subgroup for any fixed prime divisor $p$ of $|G|$.

Recall that a group is said to be $p$-closed if its Sylow $p$-subgroup is normal. In this paper, as a further generalization of Shidov's theorem, Shlyk's theorem and
{\rm\cite[Theorem 1.2]{shixuli2026}}, only considering the non-nilpotent maximal subgroups of a non-solvable group of order divisible by $p$
containing the normalizer of some Sylow subgroup, we have the following result, the proof of which is given in Section~\ref{s3}.

\begin{theorem}\ \ \label{th1} Let $G$ be a non-solvable group and $p$ any fixed prime divisor of $|G|$. Then the intersection of all non-nilpotent maximal subgroups
of $G$ of order divisible by $p$ containing the normalizer of some Sylow subgroup is $p$-closed.
\end{theorem}

Note that {\rm\cite[Theorem 1.2]{shixuli2026}} is a direct consequence of Theorem~\ref{th1}.

The following two corollaries emerge by Theorem~\ref{th1}.

\begin{cor}\ \ \label{c01} Let $G$ be a non-solvable group and $p$ any fixed prime divisor of $|G|$. Then the intersection of all non-nilpotent maximal
subgroups of $G$ of order divisible by $p$ is $p$-closed.
\end{cor}

\begin{cor}\ \ \label{c02} Let $G$ be a non-solvable group and $p$ any fixed prime divisor of $|G|$. Then the intersection of all non-normal non-nilpotent
maximal subgroups of $G$ of order divisible by $p$ is $p$-closed.
\end{cor}

For Corollary~\ref{c01}, when $p=2$, we have the following result, the proof of which is given in Section~\ref{s4}.

\begin{cor}\ \ \label{c03} The intersection of all non-nilpotent maximal subgroups of $G$ of even order is equal to the Frattini subgroup of $G$.
\end{cor}

Note that Corollary~\ref{c03} provides an extension of Shidov's theorem.

\section{Two Necessary Lemmas}\label{s2}

\begin{lemma} {\rm\cite[Theorem 10.4.2]{rob}}\ \ \label{l1} If $G$
has a nilpotent maximal subgroup of odd order, then $G$ is
solvable.
\end{lemma}

\begin{lemma} {\rm\cite[Lemma 9]{wie}}\ \ \label{l2} Let $H$ be a nilpotent Hall-subgroup of a group $G$ that is not
a Sylow subgroup of $G$. If for each prime divisor $p$ of $|H|$, assume $P\in{\rm Syl}_p(H)$ we always
have $N_G(P)=H$, then there exists a normal subgroup $K$ of $G$ such that $G=KH$ and $K\cap H=1$.
\end{lemma}

\section{Proof of Theorem~\ref{th1}}\label{s3}

\demo Since $G$ is non-solvable, $G$ has non-nilpotent maximal subgroups of order divisible by $p$ containing the normalizer of some Sylow subgroups by
{\rm\cite[Corollary 1.2]{lishi2026}}. Let $N$ be the intersection of all non-nilpotent maximal subgroups of $G$ of order divisible by $p$ containing
the normalizer of some Sylow subgroup. Then $N\trianglelefteq G$ and $N<G$.

In the following we prove that $N$ is $p$-closed.

Let $G$ be a counterexample of the smallest order. If $p\nmid|N|$. Then $N$ is naturally $p$-closed, a contradiction. Thus $p\mid|N|$.

Let $P\in{\rm Syl}_p(G)$. Then $P\cap N\in{\rm Syl}_p(N)$. Assume $P_0=P\cap N$. By the assumption, $P_0\ntrianglelefteq N$, which implies that $P_0\ntrianglelefteq G$.
One has $N_G(P)\leq N_G(P_0)<G$. Then $G$ has a maximal subgroup $T$ such that $N_G(P)\leq N_G(P_0)\leq T$. By Frattini's argument, $G=N_G(P_0)N=TN$.

Note that $p\mid|T|$. If $T$ is non-nilpotent. Then $N\leq T$. It follows that $G=TN=T$, a contradiction. Thus $T$ is nilpotent.

Claim that all nilpotent maximal subgroups of $G$ are conjugate. Let $D$ be any nilpotent maximal subgroup of $G$. Since $G$ is non-solvable, $D$ has even order by
Lemma~\ref{l1}. Let $D_2\in{\rm Syl}_2(D)$. One has $D\leq N_G(D_2)$. If $D_2$ is not a Sylow 2-subgroup of $G$. Then $D<N_G(D_2)$. It follows that $N_G(D_2)=G$ by
the maximality of $D$. One has that $D/D_2$ is a nilpotent maximal subgroup of $G/D_2$ of odd order, which implies that $G/D_2$ is solvable by Lemma~\ref{l1}. Then
$G$ is solvable, a contradiction. So $D_2$ is also a Sylow 2-subgroup of $G$. Let $E$ and $F$ be any two nilpotent maximal subgroups of $G$. Assume that $E_2$ and
$F_2$ are Sylow 2-subgroups of $E$ and $F$ respectively. Then $E_2$ and $F_2$ are also Sylow 2-subgroups of $G$. One has ${E_2}^g=F_2$ for some $g\in G$. Note that
$N_G(E_2)=E$ and $N_G(F_2)=F$. It follows that $E^g=(N_G(E_2))^g=N_G({E_2}^g)=N_G(F_2)=F$.

Therefore, any nilpotent maximal subgroup of $G$ is a conjugation of $T$.

Claim that $T$ is a nilpotent Hall-subgroup of $G$. Otherwise, if $T$ is not a Hall-subgroup of $G$. Then there exists a prime divisor $q$ of $|T|$ such that the
Sylow $q$-subgroup $Q$ of $T$ is not a Sylow $q$-subgroup of $G$. Then $Q\trianglelefteq G$. Note that $q\neq p$ and $G/Q$ is non-solvable. Let $M/Q$ be the
intersection of all non-nilpotent maximal subgroups of $G$ of order divisible by $p$ containing the normalizer of some Sylow subgroup. It is obvious that
$NQ/Q\leq M/Q$. By the minimality of $G$, one has that $M/Q$ is $p$-closed. It follows that $NQ/Q$ is also $p$-closed. Then $PQ/Q\trianglelefteq M/Q$. Note
that $PQ/Q\,{\rm char}\,NQ/Q\trianglelefteq G/Q$. One has $PQ/Q\trianglelefteq G/Q$. Then $PQ\trianglelefteq G$. By Frattini's argument,
one has $G=N_G(P)PQ=N_G(P)Q=TQ=T$, a contradiction.

Therefore, $T$ is a nilpotent Hall-subgroup of $G$. Moreover, let $U$ be any Sylow subgroup of $T$. Arguing as above $Q$, one has
$N_G(U)=T$.

In the following we discuss $T$ in three cases.

Case (1): When $p>2$. Since both $2$ and $p$ divide $|T|$, $T$ is not a Sylow subgroup of $G$. By Lemm~\ref{l2}, $G$ has a normal subgroup $H$ such
that $G=H\rtimes T$. Since $2\mid|T|$, $H$ is a Hall-subgroup of $G$ of odd order. Then $H$ is solvable by Feit-Thompson theorem. It follows that $G$ is
solvable since $G/H\cong T$ is nilpotent, a contradiction.

Case (2): When $p=2$ but $T$ is not a Sylow 2-subgroup of $G$. Arguing as Case (1), we can get a contradiction.

Case (3): When $p=2$ and $T$ is a Sylow 2-subgroup of $G$, that is $T=P$. Note that $N$ is not 2-closed by the assumption. Let $r$ be any odd prime divisor
of $|N|$. Let $R_0\in{\rm Syl}_r(N)$ and $R\in{\rm Syl}_r(G)$ such that $R_0=R\cap N$.

Claim that $R_0\trianglelefteq G$. Otherwise, if $R_0\ntrianglelefteq G$. Then $N_G(R)\leq N_G(R_0)<G$. Let $L$ be a maximal subgroup of $G$ such that
$N_G(R)\leq N_G(R_0)\leq L$. By Frattini's argument, $G=N_G(R_0)N=LN$.

We discuss $L$ in two cases.

Case $(i)$: Assume $2\mid|L|$. If $L$ is non-nilpotent. Then $N\leq L$ by the hypothesis. It follows that $G=LN=L$, a contradiction. If $L$ is nilpotent.
By above argument, $L$ is a conjugation of $T$, which implies that $L$ is also a Sylow 2-subgroup of $G$, this contradicts $r\mid|L|$ and $r>2$.

Case $(ii)$: Assume $2\nmid|L|$. Since $G=TN=PN$ and $G=LN$, one has $\frac{|P||N|}{|P\cap N|}=\frac{|L||N|}{|L\cap N|}$. Note that $(\frac{|P|}{|P\cap N|}=
\frac{|L|}{|L\cap N|})=1$. It follows that $\frac{|P|}{|P\cap N|}=1$. Then $P=P\cap N\leq N$. One has $G=PN=N$, also a contradiction.

Therefore,  $R_0\trianglelefteq G$. In particular, $R_0\trianglelefteq N$. By the arbitrariness of $r$, $N$ has a normal 2-complement, that is $N$ is 2-nilpotent.
Then $N$ is solvable by Feit-Thompson theorem. It follows that $G$ is solvable since $G=TN$ and $T$ is nilpotent, a contradiction.

All above arguments show that the counterexample of the smallest order does not exist and then $N$ is $p$-closed.\hfill\qed

\section{Proof of Corollary~\ref{c03}}\label{s4}

\demo Let $K$ be the intersection of all non-nilpotent maximal subgroups of $G$ of even order. Then $K$ is 2-closed by Corollary~\ref{c01}. It is obvious
that ${\rm\Phi}$$(G)\leq K$.

If $K\nleq{\rm\Phi}$$(G)$. Then there exists a maximal subgroup $M$ of $G$ such that $K\nleq M$. One has $G=KM$. By the hypothesis, $M$ only might be
either a non-nilpotent maximal subgroup of odd order or a nilpotent maximal subgroup which are both solvable. It follows that $G=KM$ is solvable, a contradiction.

Therefore, $K\leq{\rm\Phi}$$(G)$ and so $K={\rm\Phi}$$(G)$.\hfill\qed

\end{document}